\documentclass[10pt,a4paper]{amsart}
\input amssym.def
\input amssym.tex
\title[ Hecke groups ]
{ The Commutator subgroup of the   Hecke group  $G_5$ is not congruence}

\author{ \SMALL  cheng lien  lang}
\author{ \SMALL  mong lung lang}
\begin{document}

\baselineskip=12pt

%\subjclass{11F06; Secondary 11F03}

\maketitle

\vspace{-0.3in}

\begin{abstract} Let $q\ge 3$ be an integer
and let $G_q$ be the Hecke group associated with  $q$.
We prove that the  power subgroup $G_5^5$ and the commutator subgroup $G_5'$ are not
 congruence.

\end{abstract}

\section{Introduction}
\subsection{}
 Let $q\ge 3$ be a a fixed integer. The (homogeneous) Hecke group
 $H_q$ is defined to be the maximal discrete subgroup
of  $SL_2(\Bbb R)$
 generated by  $S$ and $T$, where $\lambda _q =2$cos$\,(\pi/q)$,
 $$ S =
\left (
\begin{array}{rr}
0 & 1 \\
-1 & 0 \\
\end{array}
\right ) \,,\,\,
 T = \left (
\begin{array}{rr}
1 & \lambda  _ q\\
0 & 1 \\
\end{array}
\right ) \,.\eqno(1.1)
$$
\noindent Let $A$ be an ideal of $\Bbb Z[\lambda_q]$.
The pricnipal congruence subgroup of $H_q$ of level $A$ is defined to be
 $$H(q, A)=\{ (a_{ij})\in H_q \,:\,a_{11} - 1, a_{22} -1, a_{12}, a_{21} \in A\}.\eqno(1.2)$$
 \noindent
 Let $Z = \left < \pm I\right >$. The (inhomogeneous) Hecke group and its
 principal congruence subgroup
  are defined as $G_q = H_q/Z$
   and $ G(q,A) = H(q,A)Z/Z$.  A subgroup $K$ of $G_q$ is congruence if $G(q, A) \subseteq K$
    for some $A$. Whether subgroups of finite indices are congruence
    have been studied extensively
    (see [F], [Lu], [S]).
  %  The problem we are interested in is that whether the  commutator subgroup $G_q'$ of $G_q$
   %  is congruent.
     In the case $q=3$, it is  known that not every subgroup
      of finite index of the modular group $G_3$
  is congruence and that the commutator subgroup $G_3'$ is congruence of level 6.
   % Proof of the congruence  of $G_3'$ is elementary  and can be
 %  viewed as a consequence of the fact that $G_3/G(3,3) \cong A_4$ has a normal
 %   subgroup of index 3 (see Lemma 3.7). As $G_q/G(q,q)$
 %    may  not possesses normal
 %    subgroup of index $q$ if $q\ge 5$,
  %     one suspects that $G_q'$ may not be congruent if $q\ge 5$.
  % The main purpose
   %   of the present article is to assure such  observation by showing that
 We suspect that $q=3$ is the only case that $G_q'$ is congruence (see Discussion 5.3).
  The main purpose
     of the present article is to show that

\noindent    {\bf Proposition 5.2.} {\em The subgroups $G_5^5$ and  $G_5'$ of the Hecke
             group $G_5$ are not congruence.}

\noindent
Note that $G_q^n$ is the subgroup of $G_q$ generated by all the elements of the
 form $ x^n\in G_q$. Note also that in the case $q\ge 3$ is a prime, these two groups
  $G_q^q$ and $G_q^2$
  are special in the sense that they are the only normal torsion subgroups of $G_q$.
 \subsection{} Our proof of the above proposition is elementary  and  requires some basic facts
  about the fundamental domains of certain subgroups of $G_5$.
   The following two facts about $G_5$ are essential in our proof as well.
\begin{enumerate}
\item[(i)] If $G_5^5$ is congruence, then $G(5,5) \subseteq G_5^5$ (Lemma 5.1).
\item[(ii)] $G_5/G(5, 5) \cong E_{5^3}PSL(2,5) $ does not possess subgroups of
 index 5 (Proposition 5.2),
 \end{enumerate}
where $E_{5^3}$ is an elementary abelian 5-group of order $5^3$. Note that
 the indices of $G_q^q$ and $G_q'$ in $G_q$ are $q$ and $2q$ respectively
 (Lemmas 3.3).

\subsection{}       The rest of the article is organised as follows. In Sections 2 and 3, we study the geometric
        aspects of the Hecke group $G_q$, such study allows us to give the geometric
         invariants (index, number of elliptic elements, number of cusps, genius)
         of $G_5^5$ and $G_5'$.
         % Note that most of the non-congruence subgroups constructed
         % in the literature are given without their geometric invariants.
                Section 4 lists all the known results which is necessary
          for our study of $G_5^5$ and $G_5'$. They are mainly results on the indices
           of the principal congruence subgroups of $G_5$. Section 5 gives us the main result
            of the present article.
 The present article is part of our project on $G_q$. We have determined
  the normalisers (see [L]) and the indices  (see [LL1], [LLT2]) for some congruence
   subgroups of $G_5$. We are currently working on the
   index formula for $G(q, \pi)$, where $q\ge 7$.

\section{Geometric invariants }
In [K], Kulkarni applied a combination of geometric
 and arithmetic methods to show that one can produce
 a set of independent generators in the sense of
 Rademacher for the congruence subgroups of the modular group,
 in fact for all subgroups of finite indices.
 His method can be generalised to all subgroups of
 finite indices of the Hecke groups $G_q$, where $q$ is a prime. See [LLT1] for detail
  (Propositions 8-10 and section 3 of [LLT1]).
 In short, for each subgroup $V$ of finite index of $G_q$,
one can associate to $V$
a set of Hecke-Farey symbols (HFS)
$\{-\infty, x_0, x_1, \cdots, x_n, \infty\},$
 a special  polygon (fundamental domain) $\Phi$,
 % and
 %a set of side parings (of this fundamental polygon) which is  a set of
 %independent generators of $S$ (determined by
 % Propositions 8-10 of [LLT1]).
 % To be more precise,
 %every subgroup $G$ of finite index
 %of $G_q1$ admits a generalised Hecke Farey Symbol (gHFS)
 %$$\{ -\infty , x_0, x_1, \cdots , x_n,   \infty \}$$
 and  an additional structure
 on each consecutive pair of $x_i$'s of the three types described
 below :
$$ {x_i} _{_{_{\smile}}} \  \hspace{-.37cm} _{_{_{_{_{_{_{\circ}}}}}}}
  x_{i+1} ,\,\,
{x_i} _{_{_{\smile}}} \  \hspace{-.37cm} _{_{_{_{_{_{_{\bullet}}}}}}}
  x_{i+1} ,\,\,
{x_i} _{_{_{\smile}} }\  \hspace{-.37cm} _{_{_{_{_{_{_{a}}}}}}}  x_{i+1}.
 $$
where $a$ is a nature  number. Each nature number $a$ occurs
 exactly twice or not at all.
Similar to the modular group,
 the actual values of the $a$'s is unimportant: it is
 the pairing induced on the consecutive pairs that matters.
 \begin{enumerate}
 \item[(i)] The side pairing  $\circ$ is an elliptic element of order 2 that pairs
  the even line $(x_i, x_{i+1})$ with itself. The trace of such an element is 0.
  \item[(ii)]
  The side pairing  $\bullet$ is an elliptic element of order $q$ that pairs
  the odd line $(x_i, x_{i+1})$. The absolute value of the trace of such an element is
   $\lambda_q$.
  \item[(iii)] The two sides with the label $a$
   are paired together by an element of infinite order.
   \item[(iv)]
The special polygon associated to the HFS is a fundamental
 domain of $V$ and the side pairings
 $ I= \{ \sigma_1,  \sigma_2 \cdots ,  \sigma_m \}$
associated to the  HFS
 is a set of independent generators of $V$
 (Theorem 7, Propositions 8-10 of [LLT1]).
 \item[(v)]
The number $d$ of  special triangles
 (a special triangle is  a fundamental domain of $G_q$) of the special polygon is  the index
 of the
   subgroup.
   \item[(vi)] The set of independent generators consists of
  $r$ matrices of infinite order, where $r$ is the number of the nature number $a$'s in the
   Hecke-Farey smybols.
   \item[(vii)] The subgroup has
   $v_2$
   (the number of the circles $\circ$ in HFS) inequivalent classes of elliptic
    elements of order 2. Each class has exactly one representative in $I$.
     \item[(viii)] The subgroup has
   $v_q$
   (the number of the bullets $\bullet$ in HFS) inequivalent classes of elliptic
    elements of order $q$. Each class has exactly one representative in $I$.
    \item[(ix)]
      The Hecke-Farey symbols
     can be partitioned into $v_{\infty}$ classes under the action of the
      set of independent generators, which gives the number of cusps of the subgroup.
    \item[(x)]  The genus $g$ can be determined by the Riemann-Hurwitz formula.
      $$(q-2)d=qv_2 +2(q-1)v_q +4qg+2qv_{\infty} -4q.\eqno(2.1)$$
     \item[(xi)] The width of a cusp $x$,
      denoted by $w(x)$,  is the number of even lines
      in $\Phi$ that comes into $x$. Algebraically, it is the smallest positive integer
       $m$ such that $\pm T_q^m$ is conjugate in $G_q$ to an element of $K$ fixing $x$
       (keep in mind that a matrix is identified with its negative in $G_q$).
       The least common multiple $N$ of the cusp widths  of $V$ is called the
        geometric width of $V$.

      \end{enumerate}
 \noindent {\bf Discussion 2.1.}
 % Note that the width of each cusp $x$ can also be determined  by investigating
  %     the special polygon $\Phi$. It is the number of even lines of $\Phi$
   %    coming into $x$.
The vertices of the Hecke-Farey symbols can be obtained by
        applying Lemma 3 of [LLT1] and the side pairings in (i)-(iii) of the above can be
         obtained by Propositions 8-10 of [LLT1].

% In particular, one can show that
%  if $U$ is a subgroup of $G_q$ ($q$ odd) of
%  index 2, then the corresponding Hecke-Farey is uniquely determined, namely
%  $$\{ {-\infty }
%  _{_{_{\smile}}} \  \hspace{-.37cm} _{_{_{_{_{_{_{\bullet}}}}}}}
%   \,\,0 \,
%   _{_{_{\smile}}} \  \hspace{-.37cm} _{_{_{_{_{_{_{\bullet}}}}}}}
%  \,\infty\}.\eqno(2.2)$$
%   It follows that $G_q$ has a unique subgroup of index 2 with independent generators
%   $$ST^{-1},\,T^{-1}S.\eqno(2.3)$$

\section{Subgroups of small indices, Power subgroups}

Let $q\ge 3 $ be a prime and let $K$ be a subgroup of $G_q$. It is clear that
 if $K$ is of index 2, then the only possible Hecke-Farey symbols for $K$ is
  $\{-\infty, 0 , \infty\}$ with the set of independent generators $\{ST^{-1}, T^{-1}S\}$.
   The invariants of $K$ is given by
    $$ d=2, v_2 = 0, v_q = 2, v_{\infty} =1, g = 0.\eqno(3.1)$$
     It is not clear that $G_q$ cannot possess subgroups of indices between 3 and $q-1$
      from algebraic point of view. However, it is clear that there is no  such Hecke-Farey
       symbols. As a consequence, we have the following :

       \smallskip
       \noindent {\bf Proposition 3.1.} {\em  Let $K$ be a subgroup of $G_q$ of index at
        most $q-1$. Then $K$ is generated by the set of independent generators  $\{ST^{-1}, T^{-1}S\}$,
        where $o(ST^{-1})= o(T^{-1}S)=q$. The invariants of $K$ are
        $ d=2, v_2 = 0, v_q = 2, v_{\infty} =1, g = 0.$  Further, $[K :  K'] = q^2$.
        }

\smallskip
\noindent {\em Proof.} Since  $\{ST^{-1}, T^{-1}S\}$ is a set of independent generators and
 $o(ST^{-1})= o(T^{-1}S)=q$, one must have  $[K :  K'] = q^2$. The rest of the proposition
  is clear. \qed

        \smallskip
        \noindent {\bf Remark.}
         Note that unlike $G_q$ ($q\ge 5$), $G_3 = PSL_2 (\Bbb Z)$ does possess subgroups of all possible
        indices, which can be proved by investigation of the Hecke-Farey symbols.

\subsection {Power subgroups of $G_q$} Denoted by $G_q^n$ the subgroup of $G_q$ generated by
 all the elements of the form $x^n$, where $x \in G_q$. It is clear that $G_q^n$ is a
  characteristic subgroup of $G_q$. Since  $G_q$ is a free product of two elliptic
   elements of orders 2 and $q$ respectively, $G_q^n$ is a proper subgroup of $G_q$ if and
    only if gcd$\,(n, 2q) \ne 1$.  The following are well known.

  \smallskip
  \noindent {\bf Lemma 3.2.} {\em Let $q$ be an odd prime.
  Then $G_q^2$ is the only subgroup of $G_q$ of index $2$. $G_q^2$ is a free
   product of two elliptic elements of order $q$. In particular,  $[G_q^2 :
   [G_q^2,G_q^2] ]= q^2$.}

  \smallskip
  \noindent {\em Proof.}
   Since $\{S, ST^{-1}\}$ is a set of independent generators of $G_q$, $o(S)=2$,
    $o(ST^{-1}) = q$, one has $ST^{-1}, T^{-1}S\in G_q^2$,
   $S \notin G_q^2.$
   We may now complete the proof of the lemma  by applying Proposition 3.1.\qed

  \smallskip
  \noindent {\bf Lemma 3.3.} {\em  Let
   $q$ be an odd prime.
   Then $G_q^q$ is the only normal subgroup of $G_q$ of index $q$.
    Further,  $G_q^q$ is a free product of $q$ elliptic elements of order $2$ and
    $[G_q^q : [G_q^q,G_q^q] ]= 2^q$.
    The invariants of $G_q^q$ are given by $d=q$, $v_2= q$, $v_q=0$,
     $v_{\infty}=1$, $g =0$.
   }

  \smallskip
\noindent {\em Proof.}  It is clear that $S\in G_q^q$, $ST^{-1}\notin G_q^q$.
 Hence $G_q^q$ is a proper subgroup that contains all the elliptic elements of order 2
 ($G_q^q$ is normal). Let $K$ be the subgroup of $G_q$ with Hecke-Farey symbols
$$\{ -\infty= {x_0 }_{_{_{\smile}} }\  \hspace{-.37cm} _{_{_{_{_{_{_{\circ}}}}}}}   {x_1}
_{_{_{\smile}} }\  \hspace{-.37cm} _{_{_{_{_{_{_{\circ}}}}}}}
x_2 , \cdots
 {x_{{q-1}/2}}
 _{_{\smile}}\  \hspace{-.36cm} _{_{_{_{_{_{_{_{\circ}}}}}}}}
  x_{{q+1}/2},\cdots
 {x_{q-2}} _{_{\smile}} \  \hspace{-.36cm} _{_{_{_{_{_{_{_{\circ}}}}}}}}
 { x_{q-1}}
  _{_{\smile}} \  \hspace{-.36cm} _{_{_{_{_{_{_{_{\circ}}}}}}}}
   x_q = \infty \},$$
    where the $x_i$' are the vertices of an ideal  $q$-gon of depth 1
     (see Discussion 2.1 of Section 2). It follows that  $[G_q : K] = q$ and that
     a set of independent generators of $K$ is given by $\{g_1, g_2, \cdots, g_{q}\}$,
      where $o(g_i) = 2$ for all $i$.
      Since $G_q^q$ contains all the elliptic elements
      of order 2, we conclude that $K$ is a subgroup of $G_q^q$. An easy study of the indices
       implies that $G_q^q = K$.
         Since $G_q^q =K$ is generated by $q$ independent generators
        of order 2, $[G_q^q : [G_q^q, G_q^q]] = 2^q.$
        Let $I$ be a normal subgroup of index $q$ of $G_q$. Since $q$ is an odd prime,
         $ S= S^q \in I$. Since $I$ is normal, $I$ contains all the elliptic elements
          of order 2. Hence $G_q^q \subseteq I$.  Since they have the same index, one must
           have $I= G_q^q$. This implies that  $G_q^q$
            is the only normal subgroup of index $q$ of $G_q$.\qed

\smallskip

\noindent {\bf Example 3.4. } The side pairings associated with the Hecke-Farey symbols of $G_5^5$ is
 given by
 {\small
 $$ \left (\begin{array}{cc}
  0 & 1\\
  -1& 0
\end{array}\right ) ,
\left (\begin{array}{cc}
  \lambda  & -1\\
  \lambda +2& -\lambda
\end{array}\right ) ,
\left (\begin{array}{cc}
 2 \lambda +1 & -2\lambda -2\\
  \lambda +2& -2\lambda-1
\end{array}\right ) ,
\left (\begin{array}{cc}
 2 \lambda +1 & -\lambda -2\\
  2\lambda +2& -2\lambda-1
\end{array}\right ) ,
\left (\begin{array}{cc}
  \lambda  & -\lambda -2\\
  1& -\lambda
\end{array}\right ) .$$}

\section { Known results about  $G_5$}
Applying the main results in [LL1] and [LL2] (Section 7 of [LL1]
 and Theorem 4.1 of [LL2]), we have the following.
\begin{enumerate}
\item[(i)] $G_5/G(5,5) \cong  G(5, \lambda +2)/G(5, 5 )  \cdot   G_5/G(5, \lambda+2)\cong  E_{5^3} PSL(2,5)$, where $ G(5, \lambda +2)/G(5, 5 )\cong E_{3^5}\cong
 \Bbb Z_5 \times \Bbb Z_5 \times \Bbb Z_5$ is the elementary
 abelian group of order $5^3$ and $ G_5/G ( 5,\lambda+2) \cong PSL(2, 5) \cong A_5$.
 \item[(ii)] Let $V$ be a congruence subgroup of $G_5$. Suppose that the
  geometric level of  $V$ is $r$ where $r$ is odd (see (xi) for the definition of the
 geometric level),  then $G(5, r)\subseteq V$.
\end{enumerate}

\section{$G_5^5$ and $G_5'$ are not congruence}

It is well known that the commutator subgroup of $\Gamma =
  G_3 $ is congruence
The main purpose of this section is to show that the commutator subgroup 
   $G_5'$ of $G_5$ is not congruence.

 \smallskip

\noindent {\bf Lemma 5.1.} {\em If $G_5^5$ is congruence, then $G(5,5)\subseteq G_5^5$.
}

\smallskip
\noindent {\em Proof.}
By Lemma 3.3, the  geometric level (see (xi) for the definition of the
 geometric level) of $G_5^5$ is 5. By (ii) of Section 4,
 $G(5,5)\subseteq G_5^5$.\qed

\subsection{} The group structure of $G(5, \lambda+2)/G(5,5)$.
 Recall first that $5 =\lambda^{-2} (\lambda +2)^2$.
By Example 3 of [LLT1],
 $$a =\left (\begin{array}{cc}
  -11\lambda - 6  & 10\lambda+5 \\
  4\lambda +3  & -4\lambda -2
\end{array}\right )=T ^{-2}
\left (\begin{array}{cc}
  3 \lambda +2 &-2\lambda-3\\
  4\lambda +3  & -4\lambda -2
\end{array}\right )\in G(5, \lambda +2) -G(5,5)\eqno(5.1)$$
By (i) of Section 4,  $G(5, \lambda+2)/G(5,5)$ is elementary abelian of order $5^3$.
 It follows that
$G(5,\lambda +2)/G(5,5)$  can be generated by $\Delta = \{ a, b=SaS^{-1},
 c= J aJ^{-1}\}$ (see (5.3) for the definition of $J$). Note that
$\Delta$ modulo $G(5,5)$  is given by
$$
a\equiv I + (\lambda+2)\left (\begin{array}{cc}
  4 & 0 \\
  4  & 1
\end{array}\right ),
b\equiv I + (\lambda+2)\left (\begin{array}{cc}
  1 & 1 \\
  0  & 4
\end{array}\right ),
c\equiv  I + (\lambda+2)\left (\begin{array}{cc}
  1 & 4 \\
  0& 4
\end{array}\right ).\eqno(5.2)$$

\smallskip
\noindent {\bf Proposition 5.2.} {\em $G_5^5$  and $G_5'$  are  not congruence.}

\medskip
\noindent {\em Proof.}  Since $G_5/G_5'$
 is abelian of order 10 and $G_5^5$ is the only normal subgroup of $G_5$ of 
  index $5$ (Lemma 3.3),  $G_5' \subseteq G_5^5$. To prove our assertion, 
   it suffices
 to show that $G_5^5$ is not congruence.
Suppose that $G_5^5$ is congruence. By Lemma 5.1,
 $G(5,5) \subseteq G_5^5.$ Since $G_5/G(5, \lambda+2)\cong A_5$ has no normal
  subgroup of index 5
  and $G_5'$ has index 5 in $G_5$, $G(5,\lambda+2)$ is not a subgroup of $G_5^5$.
    This implies that $G_5^5 G(5, \lambda+2) = G_5$.  By Second Isomorphism Theorem,
     $|G(5, \lambda+2)/[G_5^5 \cap G(5, \lambda+2)]| =5$
      and $|[G_5^5 \cap G(5, \lambda+2)]/G(5,5)|=5^2$. Note that
        $ E_{5^3} A_5 \cong G_5/G(5, 5)$ acts on
      $D = [G_5^5 \cap G(5, \lambda+2)]/G(5,5) \cong \Bbb Z_5 \times \Bbb Z_5$
       by conjugation.
       Note also that $D$ is a subgroup of $\left <\Delta \right > $.
    Recall that
    $$ J = \left (\begin{array}{cc}
  0 & 1 \\
  1& 0
\end{array}\right )\in \mbox{ Aut}\,G_5.\eqno(5.3)$$

\smallskip
\noindent Since  $[G(5,5) \cap G(5,\lambda+5)]/G(5,5) = D$  is invariant under
the conjugation of
  $E_{5^3} A_5 \cong G_5/G(5, 5)$,
$D$ is invariant under the conjugation of
 $J$ and every element of $G_5$ (in particular, $S$ and $T$).
 However, one sees by direct calculation that the only subgroup of $\left <\Delta \right >$
  invariant under $J$, $S$,    and $T$ is $\left <\Delta\right >$ itself (see Appendix A).
  A contradiction.
    Hence $G_5^5$ is not congruence. \qed

\smallskip
\noindent {\bf Discussion 5.3.} A key step in the proof of
 $G_3'$ is congruence is that
 $G_3/G(3,3)\cong A_4 \cong E_4 \Bbb Z_3$ has a normal subgroup of index 3
   (see Lemma 3.7).
   This fact is no longer true if $q =  5 $ as $G_5/G(5,5)$  possesses no normal subgroups of index
   $5$.
   As this may be true for all $q \ge 5$,
   we therefore suggest that
   $G_q'$ is not congruence if $q\ge 5$.

\medskip
\section {Appendix A}

\noindent
 {\bf Lemma A1.} {\em  Let $\pi = \lambda+2$ and let
$\Delta =\{a, b,c\}$, where
$a,b,c$ are given as in $(5.2)$. Then
the only nontrivial subgroup of
$\left <\Delta\right >$  invariant under the action of $S, T$ and $J$
 is $\left <\Delta\right >$ }.

 \smallskip
 \noindent {\em Proof.}
 Since $(I+\pi U)(I+\pi V) \equiv I+\pi(U+V)$ (mod 5), multiplication
 of $(I+\pi U)(I+\pi V)$ can be transformed into addition of $U$ and $V$.
  This makes the multiplication of
  matrices $a$, $b$, and $c$ easy. Consequently, one has
$${\tiny  r
=(ac)(ab) \equiv I + \pi
 \left (\begin{array}{cc}
   0 & 0\\
   3 & 0
 \end{array}\right ),
 s= (ac)(ab)^{-1} \equiv I + \pi
 \left (\begin{array}{cc}
   0 & 3\\
   0 & 0
 \end{array}\right ),
 t=bc  \equiv I + \pi
 \left (\begin{array}{cc}
   -3 & 0\\
   0 & 3
 \end{array}\right ).} $$
\noindent It is clear that $\left <\Delta \right >= \left <a,b,c\right >
=\left <r,s,t\right >$.
  Let $A, B \in G_5$. Set $A^B = BAB^{-1}$. Direct calculation shows that
  $$r^S = s^{-1}, r^T = rs^{-1}t^2, r^J= s,
  s^S=r^{-1}, s^T= s, s^J=r,
  t^S=t^{-1},t^T=st, t^J=t^{-1}.\eqno(A1)$$
   Denoted by $M$ a nontrivial subgroup  of $\left <r,s,t\right >$
  that is invariant under the conjugation of $J$, $S$ and $T$. Let $1\ne  \sigma = r^is^jt^k \in M$.
     One sees easily that
 \begin{enumerate}

\item[(i)] If $k \not \equiv 0$ (mod 5), without loss of generality,
 we may assume that $ k =1$.
  Then $\sigma^J \sigma^S= t^{-2} \in M$. It follows that $t\in M$.
   Hence $t^T =st \in M$. Consequently, $s\in M$. This implies
  $s^S = r^{-1} \in M$.
   In summary, $r, s, t\in M$.

 \item[(ii)]
   If $k \equiv 0$ (mod 5),
    then $\sigma$ takes the form $r^is^j$. Suppose that $i\equiv 0 $ (mod 5).
     Then $ 1\ne s^j\in M$. It follows that $s\in M$. Consequently, $r = s^T\in M$.
     Hence $rs^{-1}t^2 = r^T\in M$. As a consequence, $t \in M$.
      In summary, $r, s, t\in M$. In the case $i\not\equiv0$ (mod 5),
       we may assume that $i=1$. Hence $rs^j\in M$.
        It follows that $(rs^j)^T (rs^j)^{-1} = s^{-1}t^2 \in M$.
        Consequently, $(s^{-1}t^2)^T=  st^2 \in M$.
         This implies that  $(s^{-1}t^2)(  st^2)  = t^4 \in M$.
          Hence $t\in M$. One now sees easily that
           $r, s, t \in M$.

\end{enumerate}
\noindent Hence  the only nontrivial  subgroup of $\left< \Delta\right >$
 invariant under  $J$, $S$ and $T$ is $\left <\Delta\right>$.\qed

\bigskip

%{\small
%\noindent {\tt  lang2to46@gmail.com}
%\noindent DEPARTMENT OF MATHEMATICS,\\
%NATIONAL UNIVERSITY OF SINGAPORE,\\
%SINGAPORE 119260,\\
%REPUBLIC OF SINGAPORE\\
%}

%\noindent CHENG LIEN LANG\\
\noindent Department of Mathematics, I-Shou  University, Kaohsiung, Taiwan,
Republic of China.

\noindent   \texttt{cllang@isu.edu.tw}

\smallskip
\noindent Singapore 669608, Republic of Singapore.

\noindent \texttt{lang2to46@gmail.com}

\medskip

%\noindent {\small comm-1-4.tex}

\end{document}